\documentstyle[12pt]{article}
\begin{document}
\title{The Cardinality of the Second Uniform Indiscernible}
\author{Greg Hjorth\\      Group in Logic, University of California,
Berkeley, CA94720}
\maketitle

\begin{abstract}  
     When the second uniform indiscernible is $\aleph_{2}$, the
Martin-Solovay tree only constructs countably many reals; this
resolves a number of open questions in descriptive set theory.
\footnote {1991{\em Mathematics Subject Classification.} Primary
04A15. {\em Key words and phrases.} Descriptive set theory, uniform
indiscernibles, Martin-Solovay tree. }\\
\end{abstract}

     {\bf 1.Introduction and Definitions.}

     From now on we work in the theory ZFC + $\forall
x\epsilon\omega^{\omega}(x^{\sharp}\exists)$.

     Assuming $\forall x\epsilon\omega^{\omega}(x^{\sharp}\exists)$,
Martin and Solovay showed that every $\Sigma_{3}^{1}$ set is the
projection of a simply definable tree, T$_{2}$. This extended several
earlier results, but in some instances only with the further
assumption that T$_{2}$ has size less than $\aleph_{2}$.

     The present paper smooths the way for a further analysis by
showing that if $u_{2}$=$\aleph_{2}$ (a necessary and sufficient
condition for T$_{2}$ to have size bigger than $\aleph_{1}$) then the
smallest inner model of set theory containing T$_{2}$ and all the
ordinals, $\bf L$[T$_{2}$], has only countably many reals.

     1.1.Definition: $\gamma \epsilon${\bf Ord} is a {\em uniform
indiscernible} if $\forall x\epsilon \omega^{\omega}$($\gamma$ is an
$\bf L$[$x$] indiscernible).

     1.2.Notation: For $\alpha\epsilon${\bf Ord}, let $u_{\alpha}$ be
the $\alpha^{th}$ uniform indiscernible, beginning with
$u_{1}$=$\aleph_{1}$. For $x\epsilon\omega^{\omega}$,
$\tau\epsilon$$\cal L$($\bf L$[$x$]) indicates that $\tau$ is a skolem
function over $\bf L$[$x$], definable from no parameters.

     1.3.Representation Lemma (Solovay): For all $\gamma\epsilon$Ord,
$\beta < u_{\gamma}$, there exists $x\epsilon\omega^{\omega}$,
$\tau\epsilon$$\cal L$({\bf L}[$x$]), and $\gamma_{1}, \ldots
,\gamma_{n} < \gamma$ such that $\tau(u_{\gamma_{1}}$, $\ldots$
,$u_{\gamma_{n}}$)=$\beta$.$\Box$

     A similar result holds for finite strings of ordinals; any finite
string of ordinals less than $u_{\omega}$ can be coded from a single
real and finitely many of the uniform indiscernibles less than
$u_{\omega}$. These lemmas are signifigant because they enable us to
phrase questions about what occurs inside {\bf L}[$x$] at $\alpha <$
$u_{\omega}$ in a $\Delta_{3}^{1}$($x,y$) manner, for any
$y\epsilon\omega^{\omega}$ from which $\alpha$ can be defined using
the uniform indiscernibles.

     For our purposes it will be unimportant how the Martin-Solovay
tree is defined. It can, for example, be extracted from the scale
discussed in Moschovakis [3]. A closely related tree arises from the
``shift maps'', as implicit in the construction in 2..1. below. Both
derive from [2]. The only properties of the tree required are that:\\
\indent (i) it is no worse than $\Sigma_{3}^{1}$ in the codes for the
ordinals less than $u_{\omega}$;\\ \indent (ii) it projects to the
complete $\Sigma_{3}^{1}$ set;\\ \indent (iii) it has the same
cardinality as $u_{2}$.\\ Let us fix such a tree and call it T$_{2}$.

     1.4.Notation: For $\alpha\epsilon${\bf Ord},
$x\epsilon\omega^{\omega}$, set Next($\alpha$, $x$) to be the least
{\bf L}[$x$] indiscernible above $\alpha$.\\

{\bf 2.$u_{2}$=$\aleph_{2}$.}

2.1.Theorem(ZFC+$\forall
x\epsilon\omega^{\omega}(x^{\sharp}\exists$)): If
$u_{2}$=$\aleph_{2}$, then there are only countably many reals in $\bf
L$[T$_{2}$].

Proof:

Suppose otherwise, and now we will derive a contradiction. Let
$\theta$ be a big, regular cardinal, so that {\bf
V}$_{\theta}\models$``T$_{2}$ constructs $\geq \aleph_{1}$ many
reals''. Now choose\\
\indent (i) $N_{0}\prec N_{1}\prec \ldots N_{i} \ldots\prec${\bf V}$_{\theta}$, such that $\forall i\epsilon\omega$, $\aleph_{1}\subset N_{i}$, $\mid\mid N_{i}\mid\mid$=$\aleph_{1}$;\\
\indent (ii) ($y_{i}$)$_{i\epsilon\omega}\subset\omega^{\omega}$, $y_{i}\epsilon N_{i+1}$, $y_{i}^{\sharp}$ recursive in $y_{i+1}$, such that Next($\aleph_{1}$, $y_{i}$)$>\aleph_{2}\cap N_{i}$.\\
Let $y\epsilon\omega^{\omega}$ uniformily code the sequence
($y_{i}$)$_{i\epsilon\omega}$, let
$N_{\omega}$=$\bigcup(N_{i}$)$_{i\epsilon\omega}$, let $N$ be the
transitive collapse of $N_{\omega}$, and let T$_{2}^{*}$ be T$_{2}$ as
calculated in $N$. $N$ must be correct in calculating that T$_{2}^{*}$
constructs at least $\aleph_{1}$ many reals.  So it will suffice to
show that T$_{2}^{*}\epsilon$$\bf L$[$y$]. Let
($w_{\alpha}$)$_{\alpha\epsilon\omega_{1}}$ be a generic sequence of
reals for Coll($\omega$,$<$$\omega_{1}$).  In virtue of general facts
about forcing, it will suffice to show that T$_{2}^{*}$ is in $\bf
L$[$y$, ($w_{\alpha}$)$_{\alpha\epsilon\omega_{1}}$].  Since
(T$_{2}$)$^{N[(w_{\alpha})_{\alpha\epsilon\omega_{1}}]}$=(T$_{2}$)$^{N}$=T$_{2}^{*}$
and since every ordinal less than ($u_{\omega}$)$^{N}$ is coded by
some ($u_{1})^{N}, \ldots,$ ($u_{n})^{N}$ and ($w_{\alpha_{1}},
\ldots, w_{\alpha_{n}}, y_{0}, \ldots, y_{n}$), it suffices to show
that
Th$_{\Sigma_{3}^{1}}^{N[(w_{\alpha})_{\alpha\epsilon\omega_{1}}]}$($(w_{\alpha})_{
\alpha\epsilon\omega_{1}}\cup(y_{i})_{ i\epsilon\omega}$) can be
calculated by $\bf
L$[$y$,($w_{\alpha}$)$_{\alpha\epsilon\omega_{1}}$].

Let S be the set of $\Sigma_{3}^{1}$ sentences such that ``$\exists
z\psi(w_{\alpha_{1}}, \ldots, w_{\alpha_{n}}, y_{0}, \ldots, y_{n},
z)$''$\epsilon$S provided:\\
\indent(i) $\psi$ is $\Pi_{2}^{1}$;\\
\indent(ii) there exists $x$ recursively above $(w_{\alpha_{1}},
\ldots, w_{\alpha_{n}}, y_{0}, \ldots, y_{n})$, and a wellfounded
model $m(x)\models\exists z\psi(w_{\alpha_{1}},\ldots)\wedge${\bf
V}={\bf L}[$x$], along with indiscernibles $(c_{i})_{i\epsilon\omega ,
i>0}$ which generate $m(x)$;\\
\indent(iii) and there exists $f:${\bf
Ord}$^{m(x)}\rightarrow(u_{\omega})^{N},$ such that
$m(x)\models\gamma=\tau(c_{1}, \ldots)\epsilon [c_{m},c_{m+1})$ iff
$f(\gamma)\epsilon [(u_{m})^{N}, (u_{m+1})^{N}),$ for all
$\tau\epsilon$$\cal L$(M($x$)), and such that $f(\tau (c_{1}, \ldots,
c_{n}))=\sigma ((u_{1})^{N},
\ldots , (u_{m})^{N})$ implies that $f(\tau (c_{l_{1}}, \ldots ,
c_{l_{n}})=\sigma ((u_{l_{1}})^{N}, \ldots ,(u_{l_{m}})^{N})$, for all
$\tau\epsilon$$\cal L$($m(x)$), $\sigma\epsilon$$\cal L$($\bf
L$[$r$]), $r\epsilon
(\omega^{\omega})^{N[(w_{\alpha})_{\alpha\epsilon\omega_{1}}]}$,
$m<n$, $0<l_{1}<l_{2}, \ldots , <l_{n}<\omega$.\\ Observe that
membership in S can be phrased in terms of a tree construction, with a
given $f$ and $m(x)$ corresponding to a branch. We can view the nodes
of the tree as constructing larger and larger finite initial segments
of the theory of $m(x)[(c_{i})]_{0<i<\omega}$, as in a consistency
property, along with larger and larger finite initial segments of $f$,
that witnesses well foundedness in a particularily strong form.
Observe that the tree can in fact be calculated in {\bf
L}[$y,(w_{\alpha})_{\alpha\epsilon\omega_{1}}$], since it can locate a
set of reals $X \subset (\omega^{\omega})^{N}$ which is closed under
the pairing and sharp operations and provides codes for the ordinals
below $(u_{\omega})^{N}$. In essence, we may take $X$ to be the set of
reals generated by
$(w_{\alpha})_{\alpha\epsilon\omega_{1}}\cup(y_{i})_{i\epsilon\omega}$
by the operations of pairing and taking sharps, and observe that the
proof of Solovay's representation lemma goes through for reals rest!
ricted to this set and for ordina

Claim: if ``$\exists z\psi (w_{\alpha_{1}}, \ldots )$''$\epsilon$S,
then $N[(w_{\alpha})_{\alpha\epsilon\omega_{1}}]\models\exists z\psi
(w_{\alpha_{1}}, \ldots ).$

Given a branch ($m(x), f$), we can expand $m(x)$ out along
$\omega_{1}$ many indiscernibles, inducing a model
$m(x)(c_{\alpha})_{0<\alpha <\omega_{1}}$. By Soenfield absoluteness,
it suffices to show that this model is well founded. Now we witness
wellfoundedness by canonically extending $f$ to $f_{\omega_{1}}:{\bf
Ord}^{m(x)(c_{\alpha})_{0<\alpha <\omega_{1}}}\rightarrow
(u_{\omega_{1}})^{N}$. Given $\tau (c_{\beta_{1}} \ldots
c_{\beta_{n}})$, where $\tau\epsilon{\cal L}(m(x))$, we find
$r\epsilon (\omega^{\omega})^{N}, \sigma\epsilon{\cal L}({\bf L}[r]),
m\leq n$ such that $f(\tau (c_{1} \ldots c_{n}))=\sigma ((u_{1})^{N},
\ldots (u_{m})^{N})\epsilon [(u_{m})^{N},(u_{m+1})^{N})$, and set
$f_{\omega_{1}}(\tau (c_{\beta_{1}} \ldots c_{\beta_{n}}))=
\sigma((u_{\beta_{1}})^{N} \ldots (u_{\beta_{m}})^{N})$.

Subclaim: $f_{\omega_{1}}$ is well defined.

Suppose $\tau_{1} (c_{\beta_{1}} \ldots c_{\beta_{n}}) = \tau_{2}
(c_{\beta_{k}} \ldots c_{\beta_{n}}),$ $f(\tau_{1} (c_{1} \ldots
c_{n}))= \sigma_{1} ((u_{1})^{N} \ldots (u_{n})^{N})$ and $f(\tau_{2}
(c_{1} \ldots c_{n+1-k}))= \sigma_{2} ((u_{1})^{N} \ldots
(u_{n+1-k})^{N})$ where $\sigma_{1}$ and $\sigma_{2}$ are skolem
functions definable from reals in $N$. Then since $f$ respected the
shift maps and the $(u_{i})_{0<i<\omega}^{N}$ are joint indiscernibles
for the reals used in $\sigma_{1}$ and $\sigma_{2}$, we have that
$f(\tau_{2}(c_{k} \ldots c_{n})=f(\tau_{1}(c_{1} \ldots c_{n})=
\sigma_{2}((u_{k})^{N} \ldots (u_{n})^{N}=\sigma_{1}((u_{1})^{N}
\ldots (u_{n})^{N})$, and, consequently,
$\sigma_{2}((u_{\beta_{k}})^{N} \ldots (u_{\beta_{n}})^{N}) =
\sigma_{1}((u_{\beta_{1}})^{N} \ldots (u_{\beta_{n}})^{N})$, as
required. The other cases follow in an exactly similar fashion.

Subclaim: $f_{\omega_{1}}$ is order preserving.

This is immediate given the previous claim and the requirement that
$f$ respect the order on {\bf Ord}$^{m(x)}$.

Claim: if $N[(w_{\alpha})_{\alpha\epsilon\omega_{1}}]\models\exists
z\psi (w_{\alpha_{1}}, \ldots )$ then ``$\exists z\psi
(w_{\alpha_{1}}, \ldots )$''$\epsilon$S.

This follows by considering the appropriate sharp, and the natural map
into $(u_{\omega})^{N}$.

But since S can be calculated in {\bf L}[$y,
(w_{\alpha})_{\alpha\epsilon\omega_{1}}$], these last two claims
suffice.$\Box$

2.2.Corollaries (ZFC+$\forall
x\epsilon\omega^{\omega}(x^{\sharp}\exists)$):

(i) every $\Sigma_{3}^{1}$ set of size greater than $\aleph_{1}$
contains a perfect set;

(ii) Martin's Axiom implies $\Sigma_{3}^{1}$ lebesgue measurability;

(iii) if there is a $\Sigma_{3}^{1}$ well ordering of the reals, then
the continuum hypothesis holds;

(iv)if {\bf K}$_{0}\subset${\bf V} is an inner model with a good
$\Sigma_{3}^{1}$ well ordering of the reals, and if
$u_{2}=\aleph_{2}$, and {\bf K}$_{0}\prec_{\Sigma_{3}^{1}}${\bf V}, or
the goodness of the well ordering is in some sense sufficiently
robust, then {\bf K}$_{0}$ does not calculate $\omega_{1}$ correctly.

Proof:

These results are all proved using standard techniques along with the
fact that T$_{2}$ never constructs more than $\aleph_{1}$ many reals:
if $u_{2}=\aleph_{2}$, then this follows by the theorem; if
$u_{2}<\aleph_{2}$ then, after coding T$_{2}$ as subset of the
ordinals, this follows by the entirely general observation that a
subset of $\omega_{1}$ never constructs more than $\aleph_{1}$ many
reals.  For instance, (iii) follows as in Kechris's proof of Mansfield
theorem that if there is a $\Sigma_{2}^{1}$ well ordering of the reals
then every real is in {\bf L}.  (ii) follows by considering that
``almost every'' real must be ramdom over {\bf L}[T$_{2}$]. This
answers a question of Judah's.  $\Box$

2.3.Remark: 2.2.(iii) is rather strange sounding, but the point is
that the standard canonical inner models of large cardinals below
$\Pi_{2}^{1}$ determinacy have such ``robust'' well orders. This
indicates difficulties in forcing models of $u_{2}=\aleph_{2}$.

\end{document}